\documentclass[12pt]{report}    
\usepackage[T2A]{fontenc}
\usepackage[cp1251]{inputenc}
\usepackage[a4paper,hmargin=2.5cm,vmargin=2.5cm]{geometry}
\usepackage[english,russian]{babel}
\usepackage{graphicx}
\usepackage{placeins}
\usepackage{bmpsize}
\usepackage{amsfonts, amsmath, amssymb}
\usepackage{amsbib}

\graphicspath{}
\DeclareGraphicsExtensions{.pdf,.png,.jpg}

\def\QQ{{\mathbb Q}}
\def\ZZ{{\mathbb Z}}

\begin{document}


\begin{center}

{\Large\bfseries
Периодические замощения плоскости квадратами
\par}

\vspace{1cm}

\textbf{\large М.Д.Дмитриев}

\emph{Математический факультет Национального Исследовательского Университета Высшая Школа Экономики}

\end{center}


\textbf{Аннотация.} В работе приводится элементарное доказательство теоремы Кеньёна о том, что периодическое замощение плоскости квадратами  с периодами $(1,0)$  и $(0,\lambda)$ возможно только тогда, когда $  \lambda $ $=$ $p$ $\pm$ $\sqrt{q^2 - r^2}$ для некоторых рациональных $p \geq q \geq r \geq 0$. Доказывается аналогичный новый результат об оклеивании прямоугольника квадратами с двух сторон в один слой.

\textbf{Ключевые слова:} периодические замощения, квадрат, прямоугольник, плоскость.

\textbf{2010 MSC:} 52C20

\bigskip

Сначала сформулируем теорему, а потом введём необходимые определения.

\medskip

\textbf{ Теорема 1.} \textbf{а)}  Периодическое замощение плоскости квадратами  с периодами $(1,0)$  и $(0,\lambda)$ существует только тогда, когда $\lambda $ $=$ $p$ $\pm$ $\sqrt{q^2 - r^2}$ для некоторых рациональных $p \geq q \geq r \geq 0$.
		
		         \textbf{б)} Прямоугольник $1 \times \lambda$ можно оклеить конечным числом квадратов без просветов и наложений с двух сторон в один слой только тогда, когда  $\lambda $ $=$ $p$ $\pm$ $\sqrt{q^2 - r^2}$ для некоторых рациональных $p \geq q \geq r \geq 0$.

\medskip

\textbf{Соглашение.} Если $\lambda$ является рациональным числом, то утверждения  теоремы 1 очевидны, так как существует простой пример замощения (оклеивания). Далее мы считаем, что $\lambda$ иррационально. Также в дальнейшем будем считать, что ось $Ox$ горизонтальна, а $Oy$ вертикальна.

\medskip

\textbf{Определение.} (Ср. с \textit{[4, Определение 1]}) Зафиксируем декартову систему координат на плоскости и иррациональное число $\lambda$.

 Решетку $G  = \lbrace (m \lambda,n): m,n \in \ZZ \rbrace$ назовем \emph{порожденной} прямоугольником $1 \times \lambda $.  Введем два отношения эквивалентности.
  Назовем две точки плоскости \emph{просто эквивалентными}, если одну можно перевести в другую паралельным переносом на вектор решётки.
  Назовем две точки плоскости \emph{сложно эквивалентными}, если одну можно перевести в другую композицией центральных симметрий относительно узлов решетки $G$.

\emph{Периодическим замощением} плоскости квадратами \emph{ с периодами $(1,0)$  и $(0,\lambda)$} назовем такое конечное множество $\Gamma _1$ непересекающихся по внутренним точкам квадратов  на плоскости, что

\begin{enumerate}
\item  любая точка плоскости просто эквивалентна некоторой точке объединения квадратов из множества $\Gamma _1$;
\item    никакие две точки внутренности объединения квадратов из множества $\Gamma _1$ не просто эквивалентны.
   
\end{enumerate}
        
\emph{Оклеиванием} (с двух сторон в один слой) прямоугольника $1 \times \lambda$ квадратами назовем такое
конечное множество $\Gamma _2$ непересекающихся по внутренним точкам квадратов на плоскости,
что:
\begin{enumerate}
\item любая точка плоскости сложно эквивалентна некоторой точке объединения квадратов множества $\Gamma _2$;
\item никакие две точки внутренности объединения квадратов из множества $\Gamma _2$ не сложно эквивалентны.
\end{enumerate}

 \emph{Направлением замощения (или оклеивания)} назовём такой вектор  $u$, что стороны квадратов из множества $\Gamma_1$ (или $\Gamma _2$) либо паралельны, либо перпендикулярны $u$.

Обозначим через $2G = \lbrace (2m \lambda, 2n) : m, n \in \ZZ \rbrace $ подрешётку решётки $G$ .

Пусть дано непаралельное оси $Ox$ направление $u$ замощения (или оклеивания). Между каждыми соседними по оси $Ox$ узлами решётки $G$ (соответственно, $2G$ для теоремы 1 б) ) нарисуем \emph{ступеньку} (т.е. два перпендикулярных отрезка с общим концом, ордината которого больше ординаты этих узлов), один из отрезков которой параллелен $u$, а другой перпендикулярен (как на Рисунке 1.а). Далее проведём из каждой \emph{вершины} ступеньки (т.е. общего конца двух отрезков) отрезки, параллельные $u$, до пересечения с другой построенной ступенькой (как на Рисунке 1.б). Получим замощение плоскости $L$-образными шестиугольниками, которые могут вырождаться в прямоугольники. Назовём их \emph{уголками, построенными по вектору u для решётки $G$ (соответственно $2G$)}. Отметим, что это определение имеет смысл для любого вектора $u$, не паралельного оси $Ox$, не обязательно направления какого-то замощения или оклеивания.

Назовём уголок  \emph{разрезаемым}, если можно взять несколько квадратов (не обязательно равных), как-то разрезать каждый из них
на конечное число прямоугольников, а затем из всех получившихся прямоугольников составить данный уголок.

\includegraphics[scale=1.8]{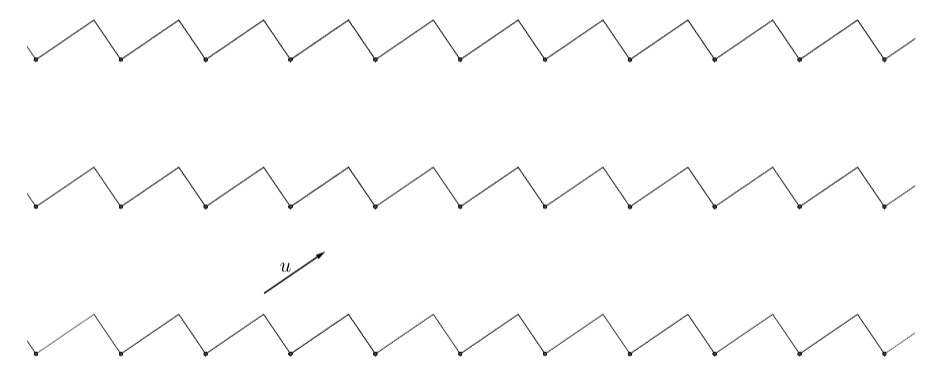}
\begin{center}
Рисунок 1.а)\textit{[3]}
\end{center}
\includegraphics[scale=2]{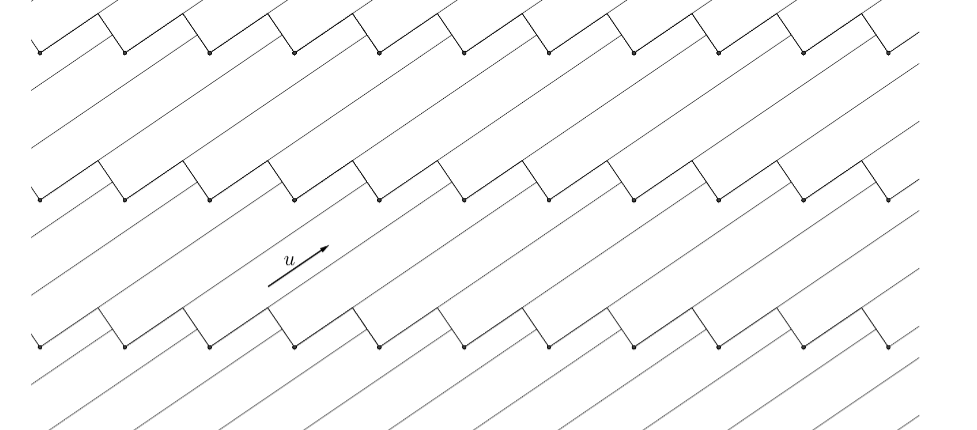}
\begin{center}
Рисунок 1.б)\textit{[3]}
\end{center}

\medskip

\medskip

\textbf{Лемма 1.} \textit{[3, Лемма 2]} \textbf{а)} Плоскость можно периодически замостить квадратами вдоль заданного направления только тогда, когда  уголок, построенный по этому направлению, разрезаем.

\textbf{б)} Если прямоугольник можно оклеить квадратами вдоль заданного направления, то уголок, построенный по этому направлению, разрезаем.

\medskip

\textbf{Доказательство леммы 1 а).}
Заметим, что уголки, которыми мы замостили плоскость, совмещаются параллельными переносами на векторы решётки $G$. Более того, \emph{линии разреза} (т.е. множество точек плоскости, эквивалентных граничным точкам квадратов из множества $\Gamma_1$ из определения замощения), проведённые внутри каждого уголка, также совмещаются. Это означает, что если квадрат на плоскости разбивается сторонами уголков на несколько частей, то равные этим частям лежат в каждом из уголков. То есть, разбив уголок на несколько частей, мы
сможем сложить из них квадраты. И так как при этом все линии разреза
будут параллельны сторонам квадратов, то и обратно - мы сможем разбить квадраты на прямоугольники и сложить из них уголок. А это и
есть определение разрезаемости. $\square$

Доказательство леммы 1 б) аналогично с заменой решётки $G$ на $2G$. 

\medskip

\textbf{Определение.} (Ср. с \textit{[3]}) Зафиксируем некоторый уголок и его разрезание на прямоугольники. Пусть уголок получается вырезанием прямоугольника $c$ $ \times $ $d$ (возможно, выраждающегося в точку) из прямоугольника $a\times b$, имеющего с первым общую вершину. Пусть  $ r_1 < r_2 < r_3 < ... < r_k$ – все длины
сторон прямоугольников разрезания. Обозначим $P = \lbrace a ; \lambda a ; d ; \lambda d ; b ; c ; $ $ r_1; \lambda r_1; r_2; \lambda r_2; r_3; \lambda r_3; ...  ; r_k; \lambda r_k\rbrace$.
Выберем специальным образом числа $t_0,t_1,t_2, ... ,t_n$  $\in P$, чтобы любое $p \in P$ единственным образом представлялось в виде $p = p_0a+p_1\lambda a+p_2t_0+...+p_{n+2} t_n$ с рациональными $p_i$ ("базис" множества $P$ над $\QQ$). Сделаем это так: выпишем в строку все числа из $P$ в указанном порядке, а потом, начиная с $\lambda a$, будем выбирать те, которые не выражаются в виде линейной комбинации с рациональными коэффициентами через предыдущие. Это и будут искомые $t_0,t_1,t_2, ...  ,t_n$. В частности, либо $t_0 = d$, либо $d = d_0a + d_1\lambda a$ с рациональными $d_0, d_1$.

  Обозначим $ P ^\prime$ $=$ $\lbrace a, \lambda a,t_0,t_1,t_2, ... ,t_n \rbrace $. Назовём число \emph{хорошим}, если оно выражается в виде линейной комбинации с рациональными коэффициентами через числа из  $ P ^\prime$. 

Пусть даны вещественное число $x$ и прямоугольник с хорошими сторонами \text{$ z_0a+z_1\lambda a+z_2t_0+...+z_{n+2}t_n$ и $ w_0a+w_1\lambda a+w_2t_0+...+w_{n+2}t_n,$} где все $z_i$ и $w_i$ рациональны. Его $x$-\emph{площадью} назовём число ($z_0+z_1x$)($w_0+w_1x$). 
Назовём $x$\emph{-площадью} уголка сумму $x$-площадей прямоугольников, на которые он разрезан. (Корректность определения показанна в работе \textit{[3, следствие из леммы 3]}.)

Следующие два свойства $x$-площади также заимствованы из работ \textit{[3]} и \textit{[5]}.

\textbf{Лемма 2.} \textit{[3, Лемма 3]}
Если прямоугольник разрезан на прямоугольники с хорошими сторонами, то его
$x$-площадь равна сумме $x$-площадей этих прямоугольников.

\medskip

\textbf{Лемма 3.} \textit{[3, Лемма 4]} Для любого действительного $x$ верно, что $x$-площадь квадрата с хорошей стороной неотрицательна при любом вещественном $x$.

\medskip

\textbf{Доказательство леммы 3.}
Если сторона квадрата равна
$z_0a+z_1 \lambda a+z_2t_0+...+z_{n+2}t_n$, где все $z_i$ - рациональны, то его $x$-площадь равна $(z_0+z_1x)^2$, что неотрицательно. $\square$

\medskip

Используем эти леммы для доказательства следующих предложений.

\medskip

\textbf{Предложение 1 а)} Если плоскость можно периодически замостить  с периодами $(1,0)$  и $(0,\lambda)$  квадратами вдоль вектора $u$, то $x$-площадь уголка, построенного по вектору  $u$, неотрицательна при любом вещественном $x$.

\textbf{б)} Если прямоугольник $1 \times \lambda$ можно оклеить конечным числом квадратов без просветов и наложений с двух сторон в один слой вдоль вектора $u$, то $x$ -площадь уголка, построенного по вектору $u$, неотрицательна при любом вещественном $x$.

\medskip

\textbf{Доказательство предложения 1 а)} По лемме 1 а), если плоскость можно замостить  с периодами $(1,0)$  и $(0,\lambda)$ вдоль вектора $u$, то разрезаем и уголок, постороенный по вектору $u$.

По лемме 2, если уголок разрезаем, то его $x$-площадь равна сумме $x$-площадей квадратов, на которые он разрезан. А так как $x$-площадь каждого квадрата неотрицательна по лемме 3, то и их сумма больше либо равна 0. Значит, $x$-площадь уголка  неотрицательна. $\square$

Аналогично с заменой леммы 1 а) на лемму 1 б) доказывается предложение 1 б).

\medskip

\textbf{Предложение 2.} Если $x$-площадь уголка, построенного в некотором направлении (для прямоугольника $1 \times \lambda$), неотрицательна при всех вещественных $x$, то $\lambda $ $=$ $p$ $\pm$ $\sqrt{q^2 - r^2}$ для некоторых рациональных $p \geq q \geq r \geq 0$. 

\medskip

Для доказательстава предложения 2 введём обозначения и докажем ещё несколько лемм.

Пусть плоскость разрезана на квадраты со сторонами, паралельными заданному направлению $u$, так, что разрезание переходит в себя при переносе на любой вектор решётки $G$. Без ограничания общности будем считать, что вектор $u$ лежит в первой координатной четверти, включая ось $Oy$, но исключая ось $Ox$ (иначе заменим его на перпендикулярный или противоположный). 
Рассмотрим уголок, посторенный по вектору $u$, вершина с наименьшей ординатой которого лежит на оси $Ox$ ближе всех к началу системы координат и отлична от него (см. Рисунок 2). Пусть он получается вырезанием прямоугольника $c \times d$ из прямоугольника $a \times b$. Введём новую прямоугольную систему координат с тем же началом, положительно ориентированную, чтобы ось ординат была параллельна вектору $u$, а ось абсцисс – перпендикулярна ему. Нетрудно проверить, что тогда узлы решётки $G$, в которых находятся вершины уголка, имеют координаты $(a,d)$ и $(c,b)$. Обозначим координаты вертикального вектора стороны прямоугольника решётки $G$ через $(e,f)$. Через 
$d_0, d_1$ и $e_0, e_1$ обозначим коэффициенты при $a$ и $\lambda a$ соответственно в представлении $d$ и $e$ в виде линейной комбинации элементов множества $P ^ \prime$ с рациональными коэффициентами.  

\includegraphics[scale=1.5]{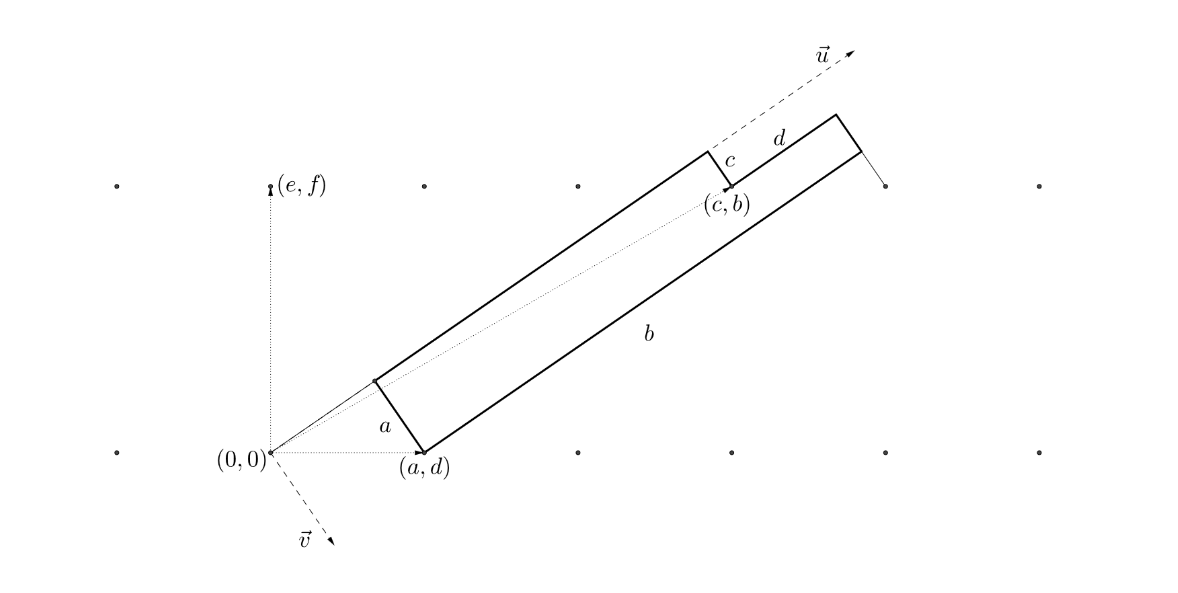}
\begin{center}
Рисунок 2.\textit{[3]}
\end{center}

 \medskip
 
\textbf{Лемма 4.} Выполнено $(e,f) = \lambda(-d,a)$ и $(c,b) = (e,f) + m(a,d)$ для некоторого целого $m$.
 
 \medskip
 
\textbf{Доказательство леммы 4.} Заметим, что вектор $(e,f)$ получается из вектора $(a,d)$ поворотом на $90^\circ$ и растяжением в $\lambda$ раз. Значит, $(e,f) = \lambda(-d,a)$. Так как вершины решётки $(a,d)$ и $(c,b)$
лежат на соседних горизонтальных уровнях решётки, то выполняется равенство между координатами соответствующих векторов
$(c,b) = (e,f) + m(a,d)$ для некоторого целого $m$.

 \medskip

\textbf{Лемма 5.} $x$-площадь уголка равна $x - (d_0 + d_1 x)(e_0 + e_1 x).$
 
 \textbf{Доказательство леммы 5.} 
Обозначим  $\widehat{z} := z_0 + z_1x$, где $z_0$ и $z_1$ - коэффициенты при $a$ и $\lambda a$ соответственно в представлении $z$ в виде линейной комбинации элементов множества $P ^ \prime$ с рациональными коэффициентами. Так как $a$ и $ \lambda a$ принадлежат $P ^ \prime$, то $ (a_0 + a_1 x)(f_0 + f_1 x) = 1 \cdot x = x$. Тогда по лемме 4 $x$-площадь уголка равна $$\widehat{a}(\widehat{b} - \widehat{d}) + \widehat{d}(\widehat{a} -\widehat{c}) = \widehat{a} \widehat{b}-\widehat{c} \widehat{d} = \widehat{a} (\widehat{f} + m \widehat{d} ) -(m \widehat{a} + \widehat{e} )\widehat{d} = \widehat{a} \widehat{f}-\widehat{e} \widehat{d} =$$  $$(a_0 + a_1 x)(f_0 + f_1 x) - (d_0 + d_1 x)(e_0 + e_1 x) = x - (d_0 + d_1 x)(e_0 + e_1 x).\square$$ 
 
 \medskip 

\textbf{Лемма 6.}  Либо $x$-площадь уголка линейно зависит от $x$, либо $d_1 \lambda ^2 + (e_1+d_0)\lambda + e_0$ $=0$ и $d_1 \not = 0$.

\medskip

\textbf{Доказательство леммы 6.} По лемме 5 $x$-площадь уголка равна $x - (d_0 + d_1 x)(e_0 + e_1 x).$

 Если $d$ или $e$ принадлежит  $P ^ \prime$, то $(d_0 + d_1 x)(e_0 + e_1 x)$ = $0$, следовательно, $x$-площадь уголка равна $x$.
 
 Если $d$ и $e$ не принадлежат  $P ^ \prime$, то $d = d_0a + d_1 \lambda a$ и $e = e_0 a + e_1 \lambda a$. Тогда по лемме 4 имеем $ - \lambda (d_0 a + d_1 \lambda a) = e_0 a + e_1 \lambda a$, следовательно,
 $ d_1\lambda ^2 + (e_1+ d_0)\lambda + e_0=0$. 

Если $d_1$ $= e_1 + d_0$ $=0$, то $e_0 =0$. Следовательно, $x$-площадь равна $x - e_1 d_0 x$ = $(1 - e_1 d_0)x$. А из того, что $e_1 + d_0$ $=0$, очевидно следует, что $1 - e_1 d_0 \not = 0$.
 
Если $d_1$ $ = 0$, а $e_1 + d_0$ $\not=0$, то $\lambda$ является рациональным числом, что противоречит соглашению после теоремы 1. $\square$

\medskip

 \textbf{Лемма 7.} Если $\lambda$ $\not=$ $p \pm \sqrt{p^2 - r^2}$ ни для каких рациональных $p$ $\geq$ $r$ $\geq 0$, то либо $x$-площадь уголка отрицательна при некотором $x$, либо выполнено неравенство $(e_0d_1 - e_1d_0 + 1)^2 - 4d_1e_0$ $< 0$.
 
\medskip
 
\textbf{Доказательство леммы 7.} По лемме 6  можно считать, что $d_1 \lambda ^2 + (e_1+d_0)\lambda + e_0$ $=0$ и $d_1$ $\not=0$. Обозначим $2p=-(e_1+d_0)/d_1$ и $q=-e_0/d_1$. Тогда $e_1= -2d_1p-d_0$ и $e_0=-d_1q$.

Тогда по лемме 5 $x$-площадь уголка равна

 $\begin{array} {rcl} S(x) & = & x - (d_0 + d_1 x)(e_0 + e_1 x) \\ & = &x - (d_0 + d_1 x)((-2d_1p-d_0)x -d_1q) \\ & = & x + (d_0 + d_1 x)(2d_1px + d_0x + d_1q) \\ & = & x^2(2d_1^2 p + d_1d_0) + x(2d_1d_0 p + d_0^2 + d_1^2q + 1) + d_1d_0q.
\end{array}$

Запишем дискриминант квадратного трёхчлена $S(x)$:

$\begin{array}{rcl} D := (2d_1d_0p + d_0^2 + d_1^2q + 1)^2-4(2d_1^2p + d_1d_0)d_1d_0q & = & (2d_1d_0p + d_0^2 - d_1^2q + 1)^2 + 4d_1^2q.\end{array}$

Рассмотрим несколько случаев:

\medskip
 
Случай 1: $D$ $> 0$. В этом случае существует  $x$, при котором $x$-площадь  отрицательна.

\medskip
 
Случай 2: $D$ $= 0$. В этом случае $-q = (2d_1d_0 p + d_0^2 - d_1^2q + 1)^2 / 4d_1^2 := r^2$ , следовательно, $\lambda$ = $p$ $\pm$ $\sqrt{p^2 + q}$ = $p$ $\pm$ $\sqrt{p^2 - r^2}$, что противоречит условию леммы. 
 
\medskip

Случай 3: $D$ $< 0$. В этом случае $q < 0$. Тогда сделаем обратную замену $2p=-(e_1+d_0)/d_1$ и $q=-e_0/d_1$. Получим, что $D = (2d_1d_0 2p + d_0^2 - d_1^2q + 1)^2 + 4d_1^2q = (e_0d_1 - e_1d_0 + 1)^2 - 4d_1e_0$ $< 0$. $\square$

\medskip

\textbf{Доказательство предложения 2.}
 Из леммы 6 мы получаем, что для некоторого направления $u$ $$d_1 \lambda ^2 + (e_1+d_0)\lambda + e_0 =0.$$ Обозначим  $e_1' = e_0d_1 / (e_1+d_0)$ и $d_0' = e_1+d_0 - e_1'$. Тогда $e_0d_1 - e_1'd_0' = (e_1')^2$ и $$(e_1+d_0)^2 - 4e_0d_1  =  (e_1'+d_0')^2 - 4e_1'd_0' - 4e_1'^2 =  (e_1'-d_0')^2 - 4e_1'^2.$$

 Cледовательно, $$\lambda  =  \frac{-(e_1+d_0)  \pm  \sqrt{(e_1+d_0)^2 - 4e_0d_1}}{2d_1} = \frac{-(e_1'+d_0') \pm \sqrt{(e_1'-d_0')^2 - 4e_1'^2}}{2d_1}   = $$ $$= -\frac{e_1'+d_0'}{2d_1} \pm \sqrt{(\frac{e_1'-d_0'}{2d_1})^2 - (\frac{e_1'}{d_1})^2} = p \pm \sqrt{q^2 - r^2},$$
   где мы обозначили $p = - \frac{e_1'+d_0'}{2d_1}$, $q=|\frac{e_1'-d_0'}{2d_1}|$ и $r=|\frac{e_1'}{d_1}|$.
   
    Докажем, что $p \geq q \geq r \geq 0$.
Неравенство $q \geq r$ очевидно следует из того, что $\lambda$ - вещественный корень квадратного уравнения.
Докажем, что $p \geq q$. Из леммы 7 следует, что $(e_0d_1 - e_1d_0 + 1)^2 - 4d_1e_0$  $< 0$. Значит, $d_1e_0 > 0$. Так как уравнение $d_1 \lambda ^2 + (e_1+d_0)\lambda + e_0$ $=0$
имеет вещественное решение, то $(e_1' + d_0')^2=(e_1+d_0)^2$ $ \geq$ $4d_1e_0$ $>0$. Заметим, что $4d_1e_0 = 4e_1' (e_1' + d_0')
$. Значит, $|e_1' + d_0'|$ $ \geq$ $4|e_1'|$. Значит, $e_1'$ и $d_0'$ одного знака и $|e_1' + d_0'|$ $ \geq$ $|e_1' - d_0'|$. Значит, $|p|$ $ \geq$ $q$, а так как $\lambda$ - положительный корень, то $p>0$. Получаем, что $p \geq q \geq r \geq 0$.
 $\square$

\medskip

Теперь мы докажем основную теорему.

 \textbf{Доказательство теоремы 1 а).} Прямо следует из предложений 1 а) и 2.
 \textbf{б).} Прямо следует из предложений 1 б) и 2. $\square$
 
 Мы предполагаем, что обратные утверждения к теореме 1 а) и б) также справедливы.


\end{document}